\documentclass[12pt,bezier]{article}
\usepackage{amssymb}
\usepackage{mathrsfs}
\usepackage{amsmath}
\usepackage{amsfonts,amsthm,amssymb}
\usepackage{amsfonts}
\usepackage{caption}
\usepackage{float}
\usepackage{graphicx}
\usepackage{graphics}
\usepackage{subfigure}
\usepackage[numbers,sort&compress]{natbib}

\newtheorem{theorem}{Theorem}
\newtheorem{lemma}{Lemma}[section]
\newtheorem*{claim}{Claim}

\begin{document}

\title{Planar graphs without cycles of length from 4 to 7 are near-bipartite}

\author{Lili Hao,\quad Weihua Yang\footnote{Corresponding author. E-mail: ywh222@163.com; yangweihua@tyut.edu.cn},\quad Shuang Zhao\\
\\ \small Department of Mathematics, Taiyuan University of
Technology,\\
\small Taiyuan Shanxi-030024,
China\\
}
\date{}
\maketitle

{\small{\flushleft \bf Abstract:}
A graph is near-bipartite if its vertex set can be partitioned into an independent set and a set which induces a forest. In this paper, planar graphs without cycles of length from 4 to 7 are shown to be near-bipartite.

\vskip 0.5cm {\flushleft \bf  Keywords}: Near-bipartite; 3-colorable; Discharging method; Planar graph

\section{Introduction}

\indent \indent A graph $G$ is \textit{$k$-degenerate} if $\delta (H)\leq k$ for each subgraph $H$ of $G$. A $k$-degenerate graph is $(k+1)$-colorable. In 1976, Borodin \cite{ref1} proposed that partitioning the vertex set of the graph into two sets whose induced graphs have better degeneracy properties. A graph $G$ is \textit{$(m,n)$-partitionable} if $V(G)$ can be partitioned into set $V_1$ and set $V_2$, where $G[V_1]$ is $m$-degenerate and $G[V_2]$ is $n$-degenerate. A graph $G$ is \textit{near-bipartite} if $V(G)$ can be partitioned into sets $V_1$ and $V_2$, where $V_1$ is independent and $G[V_2]$ is a forest. It is crystal clear that a near-bipartite graph is $(0, 1)$-partitionable. Borodin \cite{ref1} conjectured that every 5-degenerate planar graph is $(1, 2)$-partitionable and $(0, 3)$-partitionable which are confirmed by Thomassen \cite{ref10,ref11}. Clearly, a near-bipartite graph is an  extension of a 3-colorable graph in terms of shading limitations.

In 2001, Borodin and Glebov \cite{ref2} proved that every planar graph of girth at least 5 is near-bipartite, then Kawarabayashi and Thomassen \cite{ref8} made an extension of this result. Dross, Montassier and Pinlou \cite{ref6} further conjectured a stronger result that every planar graph without triangle is near-bipartite. Yang and Yuan \cite{ref22} showed that a nontrivial connected graph $G$ with $G\ncong K_4$ and $\Delta (G)\leq 3$ is near-bipartite, and it is NP-complete to judge that a graph is near-bipartite. Cranston and Yancey \cite{ref14} proved that $(1.5,-0.5)$-sparse graph without $K_4$ and Moser spindle is near-bipartite.

In 2020, Liu and Yu \cite{ref7} showed that every planar graph without $\left \{4, 6, 8 \right \}$-cycles is near-bipartite. Borodin, Glebov, Raspaud and et al \cite{ref5} showed that planar graphs without cycles of length from 4 to 7 are 3-colorable. In this paper, we proved that planar graphs without cycles of length from 4 to 7 are near-bipartite to extend above result.

The main tool used in this paper is the discharging method which has been utilized in graph theory for more than 100 years. Borodin \cite{ref15} introduced the application of discharging method in planar graph coloring, and Cranston and West \cite{ref13} explained its usage in detail.

For a class \textbf{\textit{G}} of planar graphs,  \textit{configurations} \textbf{\textit{S}}, a set of planar graphs, is \textit{unavoidable} if every graph $G$ in class \textbf{\textit{G}} has a configuration from \textbf{\textit{S}} as a subgraph.
For a given property, a \textit{reducible configuration} is a configuration which cannot occur in a minimal planar graph containing no given property. Let a type of color be the given property, i.e. solving the planar coloring problem with the above method can be translated into finding \textit{an unavoidable set of reducible configurations}. One way to prove the reducibility of configuration $C$ is to show that every coloring case of the boundary of $C$ can extend to the whole of $C$ (see Section 2). The discharging method is one way to demonstrate that a set of configurations is unavoidable, and its common  principle is as follows.
Suppose otherwise that there is a graph $G$ avoiding every configuration in \textbf{\textit{S}}, the \textit{initial charge} of $x$ is a real number associated with each vertex (and face) in $G$, such that the sum of them is negative. For some $x$, which has a negative initial charge, the charge is redistributed in such a way that all elements have non-negative charges, leading to a contradiction (see Section 3). This method has been elaborated in \cite{ref15,ref13}.

An \textit{$IF$-coloring} of a graph $G$ is a partition of its vertices into two sets $V_1$ and $V_2$ such that $V_1$, an independent set, is colored $I$ and $V_2$ colored $F$ satisfies that $G[V_2]$ is a forest. For a given graph $G$ and some cycle $C$ in $G$, an $IF$-coloring $\phi _{C}$ of $G[V(C)]$ \textit{superextends} to $G$ if there exists an $IF$-coloring $\phi _{G}$ of $G$ which entends $\phi _{C}$ satisfying there is no path joining two vertices of $C$ all of whose vertices are colored $F$. $C$ is \textit{superextendable} to $G$ if every $IF$-coloring of $G[V(C)]$ superextends to $G$. When $G$ is specify, we will say $(G, C)$ is superextendable.

A vertex $v$ is called a \textit{$k$-vertex} (\textit{$k^{+}$-vertex} or \textit{$k^{-}$-vertex}) if its degree $d(v)=k$ ($d(v)\geq k$ or $d(v)\leq k$). Denote the vertex set of neighbors of the vertex $v$ by $N(v)$. And the set of vertices locating inside (outside) the cycle $C$ is denoted by $int(C)$ ($out(C)$). A cycle $C$ is \textit{separating} if $int(C)\neq \emptyset$ and $out(C) \neq \emptyset$, otherwise $C$ is \textit{non-separating}. Let $\partial f$ denote the boundary of the face $f$.

Let $\mathcal{F}=\left \{Fb_1,Fb_2,Fb_3,Fb_4,Fb_5,Fb_6,Fb_7 \right \}$, where graph $Fb_i$ is shown in Figure 1, $i=1,2,...,7$.

\begin{figure}[htbp]
	\centering
	\subfigure[$Fb_1$]{
		\begin{minipage}[t]{0.25\linewidth}
			\centering
			\includegraphics[width=1.3in]{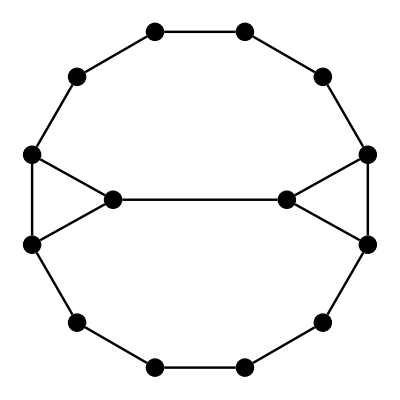}
		\end{minipage}%
		}%
	\subfigure[$Fb_2$]{
		\begin{minipage}[t]{0.25\linewidth}
			\centering
			\includegraphics[width=1.3in]{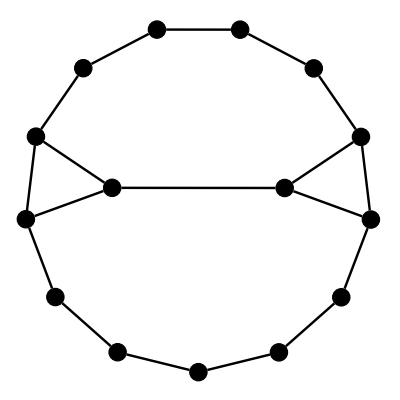}
		\end{minipage}%
		}%
	\subfigure[$Fb_3$]{
		\begin{minipage}[t]{0.25\linewidth}
			\centering							
			\includegraphics[width=1.3in]{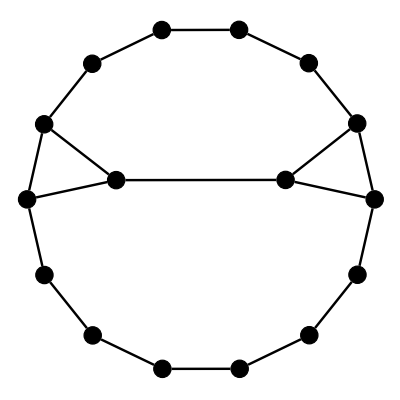}							
		\end{minipage}					
		}%
	\subfigure[$Fb_4$]{
		\begin{minipage}[t]{0.2\linewidth}
			\centering							
			\includegraphics[width=1.3in]{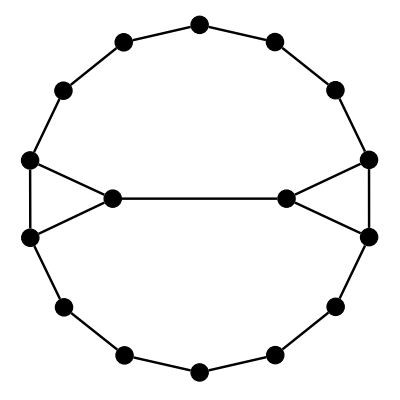}							
		\end{minipage}					
		}%

	\subfigure[$Fb_5$]{
		\begin{minipage}[t]{0.27\linewidth}
		\centering
		\includegraphics[width=1.3in]{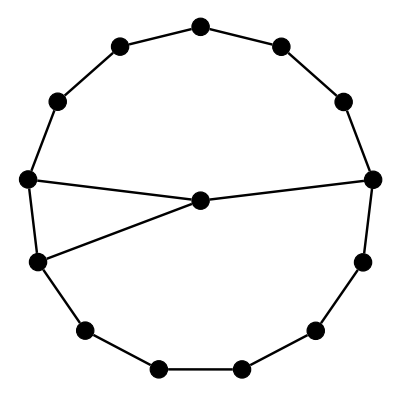}
		\end{minipage}
		}%
	\subfigure[$Fb_6$]{
		\begin{minipage}[t]{0.27\linewidth}
			\centering
			\includegraphics[width=1.3in]{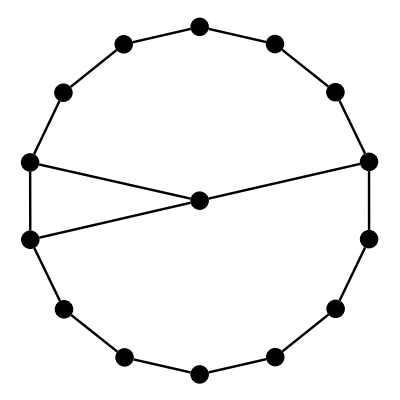}
		\end{minipage}
		}%
	\subfigure[$Fb_7$]{
		\begin{minipage}[t]{0.27\linewidth}
			\centering
			\includegraphics[width=1.3in]{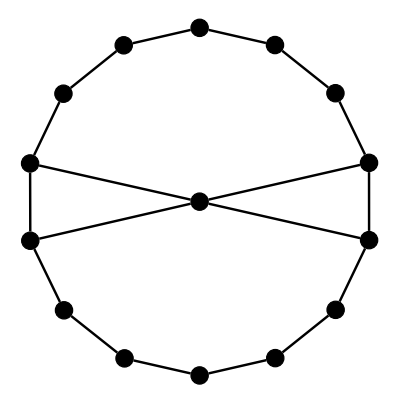}
		\end{minipage}
		}%
	
	\centering
	\caption{The graphs $Fb_i$ in $\mathcal{F}$, $i=1,2,...,7$.}

\end{figure}

We next prove the following theorem.

\begin{theorem}
	Every planar graph without cycles of length from 4 to 7 is near-bipartite.
\end{theorem}

Instead of Theorem 1, we prove a stronger result.

\begin{theorem}
	Let $G$ be a planar graph with no cycles of length from 4 to 7 and $H\notin \mathcal{F}$ for any subgraph $H$ of $G$. If $C$ is a cycle of length at most 14 in $G$, then $(G,C)$ is superextendable.
\end{theorem}

\noindent\textbf{Proof of Theorem 1. }

The proof is carried out by considering whether graph $G$ contains subgraph $Fb_i$ in $\mathcal{F}$.

\textbf{Case 1. }Suppose $Fb_1$ is a subgraph of graph $G$, then $G$ is divided into five parts $G[V(C_i)]$ with the boundary $C_i$ of the five faces of $Fb_1$ (including the exterior face), then $G[V(C_i)]=C_i$ and $(G-out(C_i),C_i)$ is superextendable according to Theorem 2, where $i=1,2,\cdots ,5$. It is evident to check that $Fb_1$ has an $IF$-coloring $\phi_{Fb_1}$. Let $\phi_{C_i}$ be equal to $\phi_{Fb_1}$ when $\phi_{Fb_1}$ is restricted to $C_i$, where $i=1,2,\cdots,5$. And then $G$ has an $IF$-coloring extended from $Fb_1$. When $Fb_j\subseteq G$, $j=2,3,\cdots ,7$, the proof is same as $Fb_1$, that is, if $H\subseteq G$ with $H\in \mathcal{F}$, then $G$ is near-bipartite.

\textbf{Case 2. }For other graphs, we can classify by girth into the following two cases. If the girth of graph $G$ is at least 8, then $G$ is near-bipartite by the result of \cite{ref2}. If the girth of $G$ is 3, then $G$ contains a 3-cycle $C$ with $G[V(C)]=C$. By Theorem 2, $C$ is superextendable to $G$, so $G$ is near-bipartite.$\hfill\blacksquare$

\section{Reducible configurations}

\indent\indent Let $(G,C_0)$ be a counterexample of Theorem 2 with minimum $\sigma(G)=|V(G)|+|E(G)|$, where $C_0$, having a precoloring, is a cycle of length at most 14 in $G$. Incontestably, $C_0$ is non-separating, otherwise, suppose that $C_0$ is separating, then $int(C_0)\neq \emptyset$ and $out(C_0)\neq \emptyset$. Therefore, $(G,C_0)$ is superextendable followed from the fact that $(G-int(C_0),C_0)$ and $(G-out(C_0),C_0)$ are superextendable, which makes a contradiction. Denote the exterior face in $G$ by $f_0$ and assume $out(C)=\emptyset$, then $C_0$ is the boundary of $f_0$. A face $f$ is \textit{internal} if $f\neq f_0$, and a vertex $v$ is \textit{internal} if $v\notin V(C_0)$, otherwise $v$ is \textit{non-internal}.

\begin{lemma}
	Let $v$ be any internal vertex in $G$, then $d(v)\geq 3$.
\end{lemma}

\noindent\textbf{Proof. }Suppose that vertex $v$ is internal with $d(v)\leq 2$ in $G$. By the minimality of $(G,C_0)$, $(G-v,C_0)$ is superextendable. If all vertices in $N(v)$ are colored $F$, then color $v$ with $I$, otherwise, color $v$ with $F$. In above cases, $C_0$ is superextendable to $G$, a contradiction.$\hfill\blacksquare$

\begin{lemma}
	There is no separating cycle of length at most 14 in $G$.
\end{lemma}

\noindent\textbf{Proof. }By contradiction. Suppose that $G$ has a separating cycle of length at most 14 denoted by $C$, then $C$ is inside of $C_0$. Firstly, by the minimality of $(G,C_0)$, $C_0$ is superextendable to $G-int(C)$, and then $G[V(C)]$ has an $IF$-coloring superextended from $C_0$. Secondly, by the minimality of $(G,C_0)$ again, $C$ is superextendable to $G-out(C)$. Thus, $C_0$ is superextendable to $G$, contrary to the hypothesis of Theorem 2.$\hfill\blacksquare$

\begin{lemma}
	$G$ is 2-connected.
\end{lemma}

\noindent\textbf{Proof. }By contradiction.

First suppose $G$ is disconnected, and $G_1$ is a component without $V(C_0)$ in $G$. By the  minimality of $(G,C_0)$, $(G-G_1,C_0)$ is superextendable. If $G_1$ contains a 3-cycle $C_1$, $(G_1,C_1)$ is superextendable by the minimality of $ (G,C_0)$. Otherwise, the girth of $G_1$ is at least 8, so $G_1$ is near-bipartite by the result of \cite{ref2}. Thus, $C_0$ is superextendable to $G$ leading to a contradiction.

Next suppose $G$ is 1-connected, and let $B$, satisfying that $B-u$ contains no vertex on $C_0$, be an endblock with a cut vertex $u$ of $G$. By the minimality of $(G,C_0)$, $(G-(B-u),C_0)$ is superextendable. Then suppose that the cycle $C_2$ is the boundary of $f$, where $f$ is the face containing $u$ with the minimum degree in $B$. If $|C_2|\leq 14$, $C_2=B[V(C_2)]$ since $C_2$ is chordless in $B$. And then by the minimality of $(G,C_0)$, $(B,C_2)$ is superextendable. If $|C_2|\geq 15$, let the neighbors of $u$ in $C_2$ be $u_1$ and $u_2$, and $C_3=uu_1u_2u$ be a 3-cycle. By the minimality of $(G,C_0)$ again, $(B+u_1u_2,C_3)$ is superextendable. In other words, whether $u$ colors $I$ or $F$, $B$ has the corresponding $IF$-coloring. Thus, $(G,C_0)$ is superextendable, a contradiction.$\hfill\blacksquare$

\begin{lemma}
	$C_0$ is chordless.
\end{lemma}

\noindent\textbf{Proof. }Assume otherwise that $C_0$ has a chord $e$. Then $G=G[V(C_0)]$, otherwise, $e$ combined a path in $C_0$ yields a separating cycle of size at most 13. Hence $C_0$ is naturally superextendable to $G$, which results in a contradiction.$\hfill\blacksquare$

Let $G\cdot \left \{u,v \right\}$ be the graph obtained by \textit{identifying} $u$ and $v$ of a graph $G$, where $u,v\in V(G)$.

An internal vertex $v$ is \textit{bad} if $v$ is incident with some 3-face and $d(v)=3$, otherwise, the internal vertex $v$ is \textit{non-bad}. An internal 4-vertex $v\in V(f)$ is \textit{poor} to the face $f$ if $v$ is incident with two 3-faces which are both adjacent to $f$, or $v$ is incident with one 3-face which is not adjacent to $f$. A vertex $v$ is \textit{willing} if $v$ is non-bad with $d(v)=3$ or poor. In a planar graph, a path $v_1v_2v_3v_4$ of the boundary of a face is a \textit{tetrad} if $v_1$, $v_2$, $v_3$, $v_4$ are bad and edges $v_1v_2$, $v_3v_4$ are incident with 3-faces as Figure 2 shows.

\begin{figure}[htbp]
	\centering
	\includegraphics[width=2.1in]{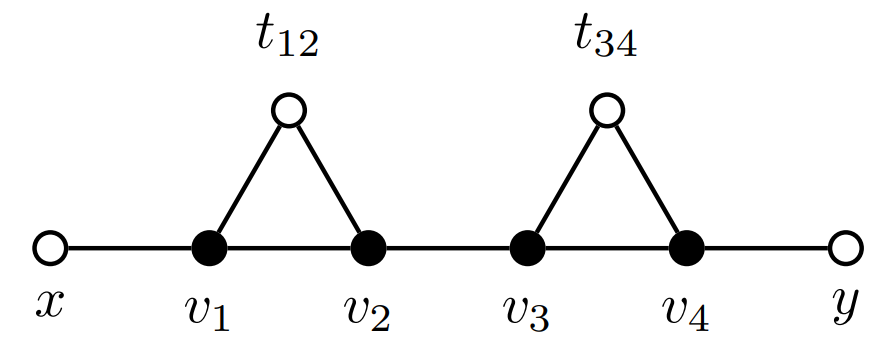}
	\caption{A tetrad}
\end{figure}

\begin{lemma}
	(1) $G$ contains no tetrad.
	
	(2) $G$ contains no face incident with five consecutive bad vertices.
	
	(3) If $G$ contains a face incident with six consecutive vertices $v_0$, $v_1$, $v_2$, $v_3$, $v_4$ and $v_5$, where $v_1$, $v_2$, $v_3$ and $v_4$ are bad, then the edges $v_0v_1$, $v_2v_3$ and $v_4v_5$ are incident with 3-faces.
\end{lemma}

\noindent\textbf{Proof. }(1) By contradiction. Suppose that $G$ contains a tetrad as Figure 2. Let $G^\ast =(G-\left \{ v_1,v_2,v_3,v_4 \right \})\cdot \left \{ x,t_{34} \right \}$. Markedly, the paths $Q,S,S_1,S_2$ mentioned below contain no edges which are absent in $G^\ast $. Suppose that all vertices of $S,S_1,S_2$ are colored $F$.

\begin{claim}
	$(G^\ast ,C_0)$ is superextendable.
\end{claim}

In fact, $G^\ast [V(C_0)]=G[V(C_0)]$, otherwise, the identifying procedure creates a chord in $C_0$ or identifies two vertices of $C_0$.  The former causes $x,{t_{34}}'\in V(C_0)$ or $t_{34},{x}'\in V(C_0)$, where ${t_{34}}'\in N(t_{34})\setminus \left\lbrace v_3,v_4\right\rbrace $ and $x'\in N(x)\setminus \left\lbrace v_1\right\rbrace $. And the latter follows that $x,t_{34}\in V(C_0)$. For above cases, there is a path of length 4 or 5, which is internally disjoint from $C_0$, between two vertices of $C_0$ denoted by $P$. Then a cycle in $P\cup C_0$ of length at most 14 separates $t_{23}$ and $t_{56}$ in $G$, contrary to Lemma 2.2. Accordingly, the precoloring of $C_0$ in $G$ remains valid in $G^\ast $.

Clearly, $G^\ast $ is a planar graph with $H\notin \mathcal{F}$ for every subgraph $H$ of $G^\ast $. $G^\ast $ contains no cycles of length from 4 to 7, otherwise, the identification procedure creates new cycles of length from 4 to 7, then there is a $(x,t_{34})$-path $Q$ of length from 4 to 7 in $G$. However, $Q\cup xv_1v_2v_3t_{34}$ is a separating cycle of length from 8 to 11, which contradicts against Lemma 2.2.

Therefore $G^\ast $, containing no subgraph $H\in \mathcal{F}$, is a planar graph without cycles of length from 4 to 7, and the claim is established by the minimality of $(G,C_0)$.

By the claim, $(G^\ast ,C_0)$ is superextendable. Next we only need to show that $G$ has an $IF$-coloring extended from any coloring of $G^\ast $ that is superextended from $C_0$, which contributes to a contradiction with hypothesis.

\textbf{Case 1.} Suppose $x$ and $t_{34}$, which are identified a new vertex in $G^\ast $, are all colored $I$. We just color $v_1,v_3,v_4$ in $F$ and $v_2$ in an opposite color of $t_{12}$.

\textbf{Case 2.} Assume $x$ and $t_{34}$ are colored $F$. Note that there is no $(x,t_{34})$-path $S$ in $G$ by the existence of an $IF$-coloring in $G^\ast $.  We color $v_2$ in $F$, $v_3$ in a same color with $y$, and $v_1$, $v_4$ in an opposite color of $t_{12}$, $y$ separately. The only condition, $t_{12}$ and $y$ are colored $F$, is invalid only when there is a $(t_{12},t_{34})$-path $S_1$. In this case, recolor $v_1$ in $F$ and $v_2$ in $I$ because of the inexistence of the $(t_{12},x)$-path $S_2$, otherwise, $S_1\cup S_2$ forms the forbidden path $S$.

(2) By contradiction. Suppose that $G$ contains a face with five consecutive bad vertices $v_1
$, $v_2$, $v_3$, $v_4$, $v_5$. By symmetry, if the edge $v_2v_3$ is incident with some 3-face, then $v_4v_5$ must be an edge on some 3-face, but $v_2v_3v_4v_5$ is a tetrad, contrary to (1). If not, there is a tetrad $v_1v_2v_3v_4$ in $G$, contrary to (1) again.

(3) Suppose that $G$ contains a face incident with six consecutive vertices $v_0$, $v_1$, $v_2$, $v_3$, $v_4$, $v_5$, where $v_1$, $v_2$, $v_3$, $v_4$ are bad. If the edge $v_2v_3$ is incident with some 3-face, then the edges $v_0v_1$, $v_4v_5$ must be on 3-faces conforming to the conclusion of Lemma 2.5(3). If not, the edges $v_1v_2$, $v_3v_4$ must be incident with 3-faces causing a tetrad $v_1v_2v_3v_4$, contrary to (1).$\hfill\blacksquare$

An 8-face $f$ with $\partial f=v_1v_2\cdots v_8$ is an \textit{$M$-face} if $v_1$, $v_2$, $v_3$, $v_5$, $v_6$, $v_7$ are bad, $v_4$ is not bad, $v_8$ is a 4-vertex, and 3-faces $t_{ij}v_iv_j$ are adjacent to $f$, where $t_{ij}=N(v_i)\cap N(v_j)$, $(i,j)\in \left\lbrace (1,8),(2,3),(5,6),(7,8)\right\rbrace $ as Figure 3(a) shows.

An 8-face $f$ with $\partial f=v_1v_2\cdots v_8$ is an \textit{$MM$-face} if $v_1$, $v_2$, $v_3$, $v_4$, $v_6$, $v_7$ are bad, $v_5$, $v_8$ are $4^{+}$-vertices, and 3-faces $t_{ij}v_iv_j$ are adjacent to $f$, where $t_{ij}=N(v_i)\cap N(v_j)$, $(i,j)\in \left\lbrace (1,8),(2,3),(4,5),(5,6),(7,8)\right\rbrace $ as Figure 3(b) shows.

\begin{figure}[htbp]
	\centering
	\subfigure[An $M$-face]{
		\begin{minipage}[t]{0.5\linewidth}
			\centering
			\includegraphics[width=2.2in]{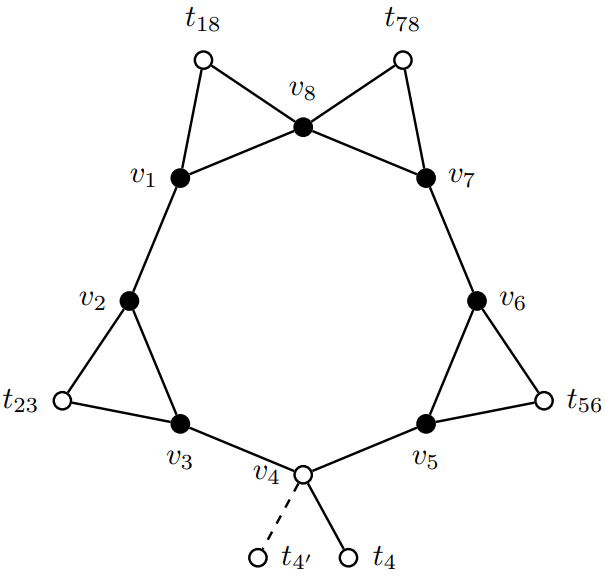}
		\end{minipage}%
	}%
	\subfigure[An $MM$-face]{
		\begin{minipage}[t]{0.5\linewidth}
			\centering
			\includegraphics[width=2.2in]{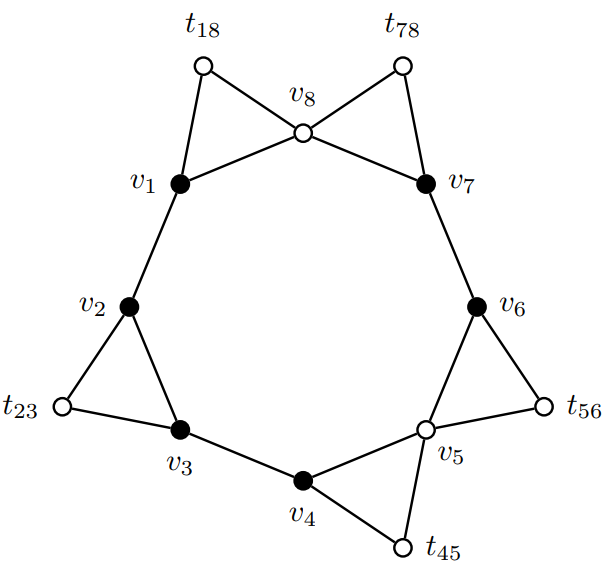}
		\end{minipage}%
	}%
	\centering
	\caption{$M$-face and $MM$-face}
	
\end{figure}

\begin{lemma}
	$G$ contains no $M$-face.
\end{lemma}

\noindent\textbf{Proof. }By contradiction. Suppose that $G$ contains an $M$-face as Figure 3(a). Let $G^\ast = (G-\left \{ v_1,v_2,v_3,v_5,v_6,v_7 \right \}+e)\cdot \left \{v_4,v_8 \right\} $, where $e=t_{18}t_{78}$. Notablel, the paths $Q_1$, $Q_2$, $S_1$, $S_2$ mentioned below contain no edges that are absent in $G^\ast $, and assume that the color is $F$ in all vertices of $S_1$, $S_2$.

\begin{claim}
	$(G^\ast ,C_0)$ is superextendable.
\end{claim}

Indeed, $G^\ast [V(C_0)]=G[V(C_0)]$, otherwise, the identification creates a chord in $C_0$ or identifies two vertices of $C_0$, or the edge $e$ is a new chord of $C_0$. The first case causes $v_4,{v_{8}}'\in V(C_0)$ or $v_8,{v_{4}}'\in V(C_0)$, where ${v_{4}}'\in N(v_4)\setminus \left\lbrace v_3,v_5\right\rbrace $ and ${v_{8}}'\in \left\lbrace t_{18},t_{78}\right\rbrace $. And in second case, it is established that $v_4,v_8\in V(C_0)$. For both cases, there exists a path $P_1$ of length 4 or 5, which is internally disjoint from $C_0$, between two vertices of $C_0$. And since $|C_0|\leq 14$, there exists a cycle contained in $P_1\cup C_0$ of length at most 14 that separates $t_{23}$ and $t_{56}$, contrary to Lemma 2.2. Third, if the adding edge $e$ is a new chord of $C_0$ in $G^\ast$, then $t_{18},t_{78}\in V(C_0)$. There is a $(t_{18},t_{78})$-path $P_2$ of length at least 6 contained in $C_0$ since the graph $G$ without cycles of length from 4 to 7. Therefore, $(C_0-P_2)\cup t_{18}v_8t_{78}$ causes a separating cycle of length at most 10, contrary to Lemma 2.2 again. So the precoloring of $C_0$ in $G$ remains valid in $G^\ast $.

Apparently, $G^\ast $ is a planar graph without subgraph $H\in \mathcal{F}$. Furthermore, $G^\ast $ contains no cycles of length from 4 to 7, otherwise, the identification or the edge-adding procedure creates new cycles of length from 4 to 7 in $G^\ast $. If the former happens, there exists a $(v_4,v_8)$-path $Q_1$ of length at most 7, so that $Q_1\cup v_8v_1v_2v_3v_4$ leads to a separating cycle of length at most 11 in $G$, contrary to Lemma 2.2. If the latter holds, there is a $(t_{18},t_{78})$-path $Q_2$ of length from 3 to 6 in $G$. Therefore, $Q_2\cup t_{18}v_1v_2v_3v_4v_5v_6v_7t_{78}$ has combined a separating cycle of size from 11 to 14 in $G$, contrary to Lemma 2.2 again.

Consequently, $G^\ast $, which contains no subgraph belonging to $\mathcal{F}$, is a planar graph without cycles of length from 4 to 7, then by the minimality of $(G,C_0)$, the claim is established .

Next we show that $G$ has an $IF$-coloring extended from the coloring of $G^\ast $. Then by the claim, $(G,C_0)$ is superextendable, a contradiction with hypothesis.

\textbf{Case 1.} Suppose that $v_4$ and $v_8$, which are identified a new vertex in $G^\ast $, are colored $I$. Afterwards, we naturally color $v_1,v_3,v_5$, $v_7$ with $F$ and $v_2$, $v_6$ with a color different from $t_{23}$, $t_{56}$ singly.

\textbf{Case 2.} Assume that $v_4$ and $v_8$ are colored $F$, then one of $t_{18}$ and $t_{78}$ is colored $I$, and the other is colored $F$ for a 3-cycle $v_8t_{18}t_{78}v_8$ in $G^\ast $. Without loss of generality, suppose that $t_{18}$ and $t_{78}$ are colored $I$ and $F$ respectively. We can then color $v_1$, $v_6$ with $F$, $v_7$ with $I$ and $v_5$ with a different color from $t_{56}$. If $t_{23}$ is colored $I$, then color $v_2$, $v_3$ with $F$. Otherwise, as $(v_4,t_{23})$-path $S_1$ and $(v_8,t_{23})$-path $S_2$ cannot coexist, one of coloring $v_2$ with $I$, $v_3$ with $F$ and the opposing coloring situation is proper.$\hfill\blacksquare$

\begin{lemma}
	$G$ contains no $MM$-face.
\end{lemma}

\noindent\textbf{Proof. }Suppose otherwise that $G$ contains an $MM$-face as Figure 3(b). Let $G^\ast = (G-\left \{ v_1,v_2,v_3,v_4,v_6,v_7 \right \}+\left \{e_1,e_2 \right \})\cdot \left \{v_5,v_8 \right\} $, where $e_1=t_{18}t_{45}$ and $e_2=t_{78}t_{56}$. And let $u$ be the vertex after identifying $v_5$ and $v_8$. It is noticeable that the paths $S_1$, $S_2$, $T$ mentioned below contain no edges which are absent in $G^\ast $ and suppose all vertices of them are colored $F$.

As in the claims of pervious two lemmas, $(G^\ast ,C_0)$ is superextendable. Next, we only need to show $G$ has an $IF$-coloring extended from the coloring of $G^\ast $, a contradiction.

\textbf{Case 1.} Assume that $v_5$ and $v_8$ are colored $I$. Then, we could color $v_1$, $v_4$, $v_6$, $v_7$ with $F$. Furthermore, since $(t_{18},t_{23})$-path $S_1$ and $(t_{23},t_{45})$-path $S_2$ cannot coexist, one of coloring $v_2$ with $F$, $v_3$ with $I$ and the opposing coloring situation is proper.

\textbf{Case 2.} Assume that $v_5$ and $v_8$ are colored $F$. As $t_{18}ut_{45}\cup e_1$ is a 3-cycle in $G^\ast $, one of $t_{18}$ and $t_{45}$ is colored $I$, and the other is $F$. The same situation occurs in $t_{56}$ and $t_{78}$. By symmetry, suppose that $t_{18}$ is colored $I$ and $t_{45}$ is colored $F$, then $(t_{45},t_{56})$-path $T$ cannot exist. So we then color $v_1$ and $v_3$ in $F$, $v_4$ in $I$, $v_2$ in a different color from $t_{23}$, and $v_6$, $v_7$ in the same color with $t_{78}$, $t_{56}$ separately.$\hfill\blacksquare$

\section{Discharging procedure}

\indent\indent The rest of proof to Theorem 2 will describe a discharging procedure, following the principles mentioned in section 1, to show an incompatible result. In this paper, we adjust the size of some initial charges such that the sum of the charges of all elements is zero, then redistribute the charges (via discharging rules) so as to the sum of the new charges of all elements is positive, which leads to a contradiction to complete the proof of Theorem 2. The specific initial charge distribution is as follows. Set the initial charge of every vertex $v\in V(G)$ with $ch(v)=d(v)-4$, any face $f\in F(G)\setminus f_0$ with $ch(f)= d(f)-4$, and $f_0$ with $ch(f_0)= d(f_0)+4$. By Euler's Formula $|V(G)|-|E(G)|+|F(G)|=2$, we have $\sum_{x\in V(G)\cup F(G)}ch(x)=0$.

Let $v\in V(C_0)$ be a 3-vertex, and $f_1$, $f_2$ be its incident faces except $f_0$. The vertex $v$ is \textit{content} if $f_1$ is a 3-face and $f_2$ is an $8^{+}$-face.

An 8-face $f$ with $\partial f=v_1v_2\cdots v_8$ is an \textit{$Fa_1$-face} if $v_1,v_3$ are content, $v_2$ is a 2-vertex, $v_4,v_5,v_6,v_8$ are bad, $v_7$ is willing, and 3-faces $t_{ij}v_iv_j$ are all adjacent to $f$, where $t_{ij}=N(v_i)\cap N(v_j)$, $(i,j)\in \left\lbrace (1,8),(3,4),(5,6)\right\rbrace $ and $t_{34},t_{18}\in V(C_0)$ as Figure 4(a) shows.

An 8-face $f$ with $\partial f=v_1v_2\cdots v_8$ is an \textit{$Fa_2$-face} if $v_1,v_3$ are content, $v_2$ is a 2-vertex, $v_4,v_5,v_7,v_8$ are bad, $v_6$ is willing and 3-faces $t_{ij}v_iv_j$ are adjacent to $f$, where $t_{ij}=N(v_i)\cap N(v_j)$, $(i,j)\in \left\lbrace (1,8),(3,4),(5,6),(6,7)\right\rbrace $ and $t_{34},t_{18}\in V(C_0)$ as Figure 4(b) shows.

\begin{figure}[htbp]
	\centering
	\subfigure[$FA_1$]{
		\begin{minipage}[t]{0.5\linewidth}
			\centering
			\includegraphics[width=2in]{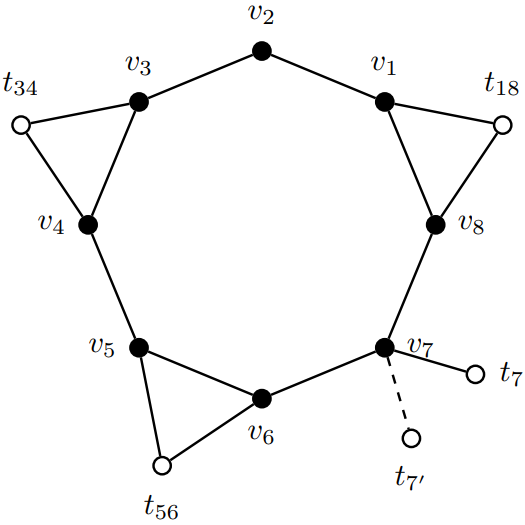}
		\end{minipage}%
	}%
	\subfigure[$FA_2$]{
		\begin{minipage}[t]{0.5\linewidth}
			\centering
			\includegraphics[width=2in]{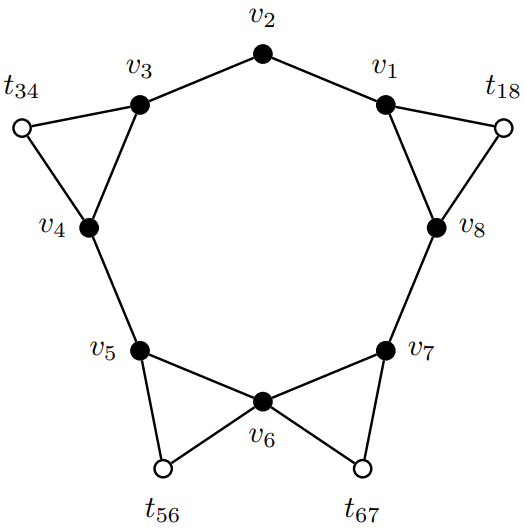}
		\end{minipage}%
	}%
	\centering
	\caption{$FA$-structures}
	
\end{figure}

Let $FA_i$ denote the subgraph $G[V(f)\cup N(V(f))]$ of $G$, where $f$ is an $Fa_i$-face, $i=1,2$. We call $FA_1$ or $FA_2$ \textit{$FA$-structure}. The vertex $v\in V(f)$ is \textit{special} to the face $f$ if $v$ is content and $f$ is an $Fa_1$-face or $Fa_2$-face.

The discharging procedure meets the following discharging rules:

\noindent \textbf{R1.} Each 3-face receives $\frac{1}{3}$ from each incident vertex.

\noindent \textbf{R2.} Each internal $8^{+}$-face $f$ sends to each incident vertex $v$:

\textbf{R2.1.} $\frac{2}{3}$ if $d(v)=2$ or $v$ is bad;

\textbf{R2.2.} $\frac{1}{3}$ if $v$ is willing or content but not special.

\noindent\textbf{R3.} Each internal $8^{+}$-face $f$ receives $\frac{1}{3}$ from each incident vertex $v$ if:

\textbf{R3.1.} $v$ is an interal $5^{+}$-vertex incident with two 3-faces adjacent to $f$;

\textbf{R3.2.} $v$ is a non-internal $4^{+}$-vertex.

\noindent \textbf{R4.} The exterior face $f_0$ sends to each incident vertex $v$:

\textbf{R4.1.} $\frac{4}{3}$ if $d(v)=2$ or $v$ is special;

\textbf{R4.2.} $1$ if $v$ is not special with $d(v)\geq 3$.

\noindent \textbf{R5.} Each internal face sends the surplus charge to $f_0$.

Let $ch^{\ast }(x)$ and  $ch_1(x)$ be the charge of $x\in V(G)\cup F(G)$ after the discharging procedure (i.e. \textit{final charge}) and after the discharging procedure except R5, respectively.

\begin{lemma}
	$ch^{\ast }(v)\geq 0$ when $v\in V(G)$.
\end{lemma}

\noindent\textbf{Proof. }The case is divided by considering whether $v$ is internal.

\textbf{Case 1.} Suppose that $v\in V(C_0)$.

If $d(v)=2$, then its incident internal face $f$ has $d(f)\geq 8$, otherwise, $(G,C_0)$ is superextendable as $G=G[V(C_0)]$. By R4.1 and R2.1, $v$ gets $\frac{4}{3}$ from $f_0$ and $\frac{2}{3}$ from $f$, and thence $ch^{\ast}(v)=2-4+\frac{2}{3}+\frac{4}{3}=0$.

If $d(v)=3$, $v$ is incident to $f_0$ and two internal faces $f_1$ and $f_2$. First, suppose that $v$ is content, let $f_1$ be a 3-face. By R1, R2 and R4, if $v$ is special to $f_2$, $v$ sends $\frac{1}{3}$ to $f_1$ and gets $\frac{4}{3}$ from $f_0$, or sends $\frac{1}{3}$ to $f_1$ and gets $\frac{1}{3}$ from $f_2$ and 1 from $f_0$ otherwise. Next, considering that both $f_1$ and $f_2$ are $8^{+}$-faces, then $v$ only gets 1 from $f_0$. Hence, $ch^{\ast}(v)\geq 3-4+\min\left\lbrace -\frac{1}{3}+\frac{4}{3},-\frac{1}{3}+\frac{1}{3}+1,1\right\rbrace =0$.

If $d(v)\geq 4$, then $v$ gets 1 from $f_0$ by R4.2 and sends $\frac{1}{3}$ to each of the $d(v)-1$ incident internal faces by R1 and R3.2. Consequently, $ch^{\ast}(v)= d(v)-4+1-\frac{1}{3}\times (d(v)-1)=\frac{2}{3}(d(v)-4)\geq 0$.

\textbf{Case 2.} Suppose that $v\notin V(C_0)$, then $d(v)\geq 3$ by Lemma 2.1. Denote that $n_3(v)$ is the number of 3-faces incident with $v$.

If $d(v)= 3$, $0\leq n_3(v)\leq 1$ since there is no 4-cycles in $G$. If $n_{3}(v)=0$, then $v$ is willing and gets $\frac{1}{3}$ from each of the three incident $8^{+}$-faces by R2.2. If $n_{3}(v)=1$, then $v$ is bad. By R1 and R2.1, $v$ is incident to one 3-face getting $\frac{1}{3}$ from $v$ and two $8^{+}$-faces sending $\frac{2}{3}$ to $v$, that implies $ch^{\ast}(v)\geq 3-4+\min \left\lbrace \frac{1}{3}\times 3,-\frac{1}{3}+\frac{2}{3}\times 2\right\rbrace = 0$.

If $d(v)= 4$, distinctly, $0\leq n_3(v)\leq 2$. Let the faces incident with $v$ be named $f_1,f_2,f_3,f_4$ in order. If $n_{3}(v)=0$, there is no charge exchange between $v$ and $f_i$, $i=1,2,3,4$. If $n_{3}(v)=1$, let $f_1$ be a 3-face, then $v$ is poor to $f_3$. So, $v$ sends $\frac{1}{3}$ to $f_1$ by R1 and gets $\frac{1}{3}$ from $f_3$ by R2.2. If $n_{3}(v)=2$, let $f_1,f_3$ be 3-faces, then $v$ is poor to $f_2$ and $f_4$. Therefore, $v$ sends $\frac{1}{3}$ to $f_1,f_3$ by R1 and gets $\frac{1}{3}$ from $f_2$, $f_4$ by R2.2. Hence, $ch^{\ast}(v)\geq 4-4+\min\left\lbrace 0,-\frac{1}{3}+\frac{1}{3},-\frac{1}{3}\times 2+\frac{1}{3}\times 2\right\rbrace = 0$.

If $d(v)= 5$, $0\leq n_{3}(v)\leq 2$. Then $v$ loses the most charge in case when $n_{3}(v)=2$, $v$ sends $\frac{1}{3}$ to each of its incident two 3-faces $f_1$, $f_2$ and one $8^{+}$-face adjacent with both $f_1$ and $f_2$ by R1 and R3.1. Thus, $ch^{\ast}(v)\geq 5-4-\frac{1}{3}\times 2-\frac{1}{3}= 0$.

If $d(v)\geq 6$, the case where $v$ gives the most charge is $\frac{1}{3}$ to all its incident faces by R1 and R3.1. Thereby $ch^{\ast}(v)\geq d(v)-4-\frac{1}{3}\times d(v)=\frac{2}{3}(d(v)-6)\geq 0$.

For all above cases, $ch^{\ast}(v)\geq 0$ when $v\in V(G)$.$\hfill\blacksquare$

\begin{lemma}
	$ch^{\ast }(f)\geq 0$ when $f\in F(G)\setminus f_0$.
\end{lemma}

\noindent\textbf{Proof. }It is obvious that we only need to show $ch_{1}(f)\geq 0$.

\textbf{Case 1.} If $d(f)=3$, then $f$ gets $\frac{1}{3}$ from each incident vertex by R1, which represents $ch_1(f)=3-4+\frac{1}{3}\times 3=0$.

\textbf{Case 2.} Suppose that $f$ is an $8^{+}$-face having some 2-vertices, in that way, $f$ contains two common $3^{+}$-vertices $u,v$ with $f_0$. Let $n_2(f)$ denote the number of 2-vertices in $f$. Observe that $1\leq n_2(f)\leq d(f)-3$, and a content vertex cannot be special when $d(f)=8$ and $n_2(f)\geq 2$ or $d(f)\geq 9$. Notably, $f$ receives $\frac{1}{3}$ from $w\in V(C_0)$ if $d(w)\geq 4$. By R2, $f$ gives the most charge to its incident vertices in cases where $\frac{2}{3}$ is given to each of the $d(f)-2$ vertices which are 2-vertices or bad, and $\frac{1}{3}$ to two content but not special vertices.
Therefore, $ch_1(f)\geq d(f)-4-\frac{2}{3}\times (d(f)-2)-\frac{1}{3}\times 2=\frac{1}{3}(d(f)-10)$. Let $\partial f=vv_1\cdots v_{n_2(f)}uu_1\cdots  u_{d(f)-n_2(f)-2}$, $d(u)=3$ with $N(u)=\left\lbrace u_1,v_{n_2(f)},b\right\rbrace $ and $d(v)=3$ with $N(v)=\left\lbrace v_1,u_{d(f)-n_2(f)-2},a \right\rbrace $, where $v_i$ is a 2-vertex, $u_j$ is internal, $a,b\in V(C_0)$, $1\leq i\leq n_2(f)$ and $1\leq j\leq d(f)-n_2(f)-2$. If $u$ (or $v$) is content, let $f_u$ (or $f_v$) be its incident 3-face with $\partial f_u=uu_1b$ (or $\partial f_v=vu_{d(f)-n_2(f)-2}a$).

\textbf{Subcase 2.1.} When $d(f)\geq 10$,  $ch_1(f)\geq 0$ and $ch_1(f)=0$ only if $d(f)=10$ and $f$ is incident with eight 2-vertices or bad vertices and two content vertices $u$, $v$.

If $n_2(f)=1,2,3$, $ch_1(f)>0$ since there is at least one non-bad vertex by Lemma 2.5(2).

If $n_2(f)=4$ and $G$ satisfies the condition of $ch_1(f)=0$, then $u_2u_3$ is an edge of a 3-face $f_1$. If $V(f_1)\cap V(C_0)=\emptyset$, there exists a cycle contained in $au_4u_3u_2u_1b\cup C_0$ of length at most 12, contrary to Lemma 2.2. If $V(f_1)\cap V(C_0)=t_{23}$, there is a separating cycle in $au_4u_3t_{23}\cup C_0$ or $bu_1u_2t_{23}\cup C_0$ of length from 4 to 7, contrary to the restrictions of Theorem 2.

If $n_2(f)=5,7$, $ch_1(f)>0$ since a non-bad vertex exists in $f$.

If $n_2(f)=6$, $V(G)=V(C_0)\cup \left \{ u_1,u_2 \right \}$, otherwise there is a separating cycle contained in $au_2u_1b\cup C_0$ of length 8, contrary to Lemma 2.2. Considering that $G$ has no $\left \{4,5,6,7 \right \}$-cycles, $G\cong Fb_3\in \mathcal{F}$, a contradiction.

\textbf{Subcase 2.2.} When $d(f)=9$, $ch_1(f)<0$ only if $f$ has seven 2-vertices or bad vertices  and two content vertices $u$, $v$.

If $n_2(f)=1,2$, $ch_1(f)\geq 0$ since there is at least one non-bad vertex in $f$ by Lemma 2.5(2).

If $n_2(f)=3$ and $G$ satisfies the condition of $ch_1(f)< 0$, then edge $u_2u_3$ is on a 3-face $f_1$. As mentioned in Subcase 2.1, there is a contradiction whether $V(f_1)\cap V(C_0)$ is empty or not.

If $n_2(f)=4,6$, $ch_1(f)\geq 0$ as there exists a non-bad vertex in $f$.

If $n_2(f)=5$, as in Subcase 2.1, $G\cong Fb_2$ or $Fb_4\in \mathcal{F}$, which makes a contradiction.

\textbf{Subcase 2.3.} If $d(f)=8$, then $ch_1(f)<0$ only if $f$ satisfies one of the following two conditions by R2. Con 1 is that $f$ has six 2-vertices or bad vertices receiving $\frac{2}{3}$ from $f$, one content but not special vertex $u$ receiving $\frac{1}{3}$ from $f$, and in this case $v$ must be a 3-vertex. Con 2 is that $f$ contains five 2-vertices or bad vertices, two content but not special vertices $u$, $v$ and one willing vertex.

When $n_2(f)=1$, discuss the existence of Con 1 and Con 2 correspondingly. Con 1 cannot be satisfied as Lemma 2.5(2). Consider Con 2, $f$ contains at least one non-bad vertex causing an $Fa_1$-face or $Fa_2$-face, so content vertices $u,v$ are special to $f$, and Con 2 cannot exist.

Assume that $n_2(f)=2$. In Con 1, if there is one content but not special vertex $u$ and four bad vertices $u_1$, $u_2$, $u_3$, $u_4$ in $f$, then $v$ must be content since $vu_4$ is on some 3-face by Lemma 2.5(3). Hence, $f$ sends $\frac{2}{3}$ to two 2-vertices and four bad vertices by R2.1, and $\frac{1}{3}$ to $u$, $v$ by R2.2. Subsequently, $ch_1(f)=8-4-\frac{2}{3}\times 6-\frac{1}{3}\times 2=-\frac{2}{3}$, but there is a separating cycle in $au_4u_3u_2u_1b\cup C_0$ of length at most 14, contrary to Lemma 2.2. Con 2 also cannot exist since a separating cycle with length at most 14 is shown as in Con 1, and $u_2$ or $u_3$ cannot be willing.

Assume that $n_2(f)=3$. Both the establishment of Con 1 and Con 2 yields a separating cycle in $vu_3u_2u_1u\cup C_0$ of length at most 14, contrary to Lemma 2.2.

Assume that $n_2(f)=4$. As $G$ contains no cycles of length from 4 to 7, Con 1 is established only when $f$ has two bad vertices $u_1,u_2$, and suppose $f$ has one content vertex $u$, then $uu_1$, $vu_2$ is in a 3-faces and $v$ must be content. Now, $G$ contains no internal vertiex besides $u$ and $v$, otherwise there is a separating cycle in $au_2u_1b\cup C_0$ of length at most 10. Therefore, $G\cong Fb_1,Fb_2$ or $Fb_3$ that brings about a contradition to the conditions of Theorem 2. Con 2 is not established since $u_1$ or $u_2$ cannot be willing when both $v$ and $u$ satisfy the condition.

Assume that $n_2(f)=5$. Con 1 needs one bad vertex and one content vertex $u$, then $V(G)=V(C_0)\cup \left\lbrace u_1\right\rbrace $, otherwise, part of $vu_1b\cup C_0$ is a separating cycle of length 8 or 9. Then $G\cong Fb_5$ or $Fb_6 \in \mathcal{F}$, contrary to conditions for Theorem 2. Con 2 holds only if both $u$ and $v$ are content and $u_1$ is poor, then as in Con 1, $G\cong Fb_7$, contrary to Theorem 2 again.

\textbf{Case 3.} Suppose that $f$ is an $8^{+}$-face that shares some $3^{+}$-vertices with $f_0$ but no 2-vertex. Let $n_{3^{+}}(f)$ denote the number of non-internal vertices in $f$ and $\partial f=v_1v_2\cdots v_{d(f)}$. Obviously $1\leq n_{3^{+}}(f) \leq 2$.

\textbf{Subcase 3.1.} If $n_{3^{+}}(f)=1$, then assume $v_1\in V(C_0)$ and $v_i\notin V(C_0)$, $2\leq i\leq d(f)$. Due to Lemma 2.5(2), $f$ gets $\frac{1}{3}$ from $v_1$ as $d(v_1)\geq 4$ and at most sends $\frac{2}{3}$ to bad vertices with cardinality $d(f)-2$, and $\frac{1}{3}$ to one willing vertex by R2 and R3.2. Hence, $ch_{1}(f)\geq d(f)-4+\frac{1}{3}-\frac{2}{3}\times (d(f)-2)-\frac{1}{3}\times 1=\frac{1}{3}(d(f)-8)\geq 0$.

\textbf{Subcase 3.2.} If $n_{3^{+}}(f)=2$, let $v_1,v_2\in V(C_0)$ and $v_i\notin V(C_0)$, $3\leq i\leq d(f)$. Owing to Lemma 2.5(2), $f$ gives at most $\frac{2}{3}$ to each of $d(f)-3$ bad vertices, $\frac{1}{3}$ to one willing vertex and $\frac{1}{3}$ to each of $v_1$, $v_2$ when they are content by R2. Notably, a content vertex is not special in this case. Therefore, $ch_{1}(f)\geq d(f)-4-\frac{2}{3}\times (d(f)-3)-\frac{1}{3}\times 1-\frac{1}{3}\times 2=\frac{1}{3}(d(f)-9)$, and $ch_{1}(f)<0$ only if $f$ is an 8-face with five bad vertices, one willing vertex and two content vertices $v_1$, $v_2$. The grouping of bad vertices can be divided into the following two types, group 1 is 1-4 and group 2 is 2-3. In group 1, suppose that all internal vertices are bad except for $v_4$ that is willing. By Lemma 2.5(3), $v_1v_8$, $v_6v_7$, $v_4v_5$ are incident with 3-faces, at this time $v_4$ violates its requirements of group 1. In group 2, assume that $v_3$, $v_4$, $v_6$, $v_7$, $v_8$ are bad and $v_5$ is willing, then $v_4v_5$, $v_6v_7$ are incident with 3-faces, but $v_5$ cannot be willing to result in the inactivity of group 2.

\textbf{Case 4.} Suppose that $d(f)\geq 8$ and $f$ shares no vertex with $f_0$. By Lemma 2.5(2), there are at least two non-adjacent non-bad vertices in $f$, thus by R2 $f$ gives to at most $d(f)-2$ incident bad vertices $\frac{2}{3}$ and at least two other incident vertices $\frac{1}{3}$ when they are willing. Accordingly, it follows that $ch_1(f)\geq d(f)-4-\frac{2}{3}\times (d(f)-2)-\frac{1}{3}\times 2=\frac{1}{3}(d(f)-10)$. Let $\partial f=v_1v_2\cdots v_{d(f)}$.

\textbf{Subcase 4.1.} If $d(f)\geq 9$, $ch_1(f)<0$ only if $f$ exactly has seven bad vertices and two willing vertices, then all bad vertices have to be grouped by 3-4. Suppose that $v_4,v_9$ are willing, while others are bad. Since Lemma 2.5(3), the edges $v_4v_5,v_6v_7$, $v_8v_9$ are incident with 3-faces. Afterwards, $v_4, v_9$ must be 4-vertices otherwise they are bad, but one of them cannot be poor, which occasions $ch_1(f)\geq 0$.

\textbf{Subcase 4.2.} Let $d(f)=8$, $ch_1(f)<0$ only if $f$ satisfies one of the following two conditions. Con 1 is that $f$ contains six bad vertices and at least one willing vertex. Con 2 is that $f$ contains five bad vertices and three willing vertices.

Con 1. In $f$, all bad vertices can be divided into two groups, group 1 is 2-4 and group 2 is 3-3. In group 1, suppose that all vertices on $f$ except $v_1,v_4$ are bad, then edges $v_4v_5$, $v_6v_7$, $v_1v_8$ are on 3-faces as Lemma 2.5(3). If the edge $v_2v_3$ is incident with some 3-face, then neither $v_1$ nor $v_4$ is willing. Otherwise, the edges $v_1v_2$, $v_3v_4$ are on 3-faces taking an $MM$-face, which makes a contradiction with Lemma 2.7. In group 2, suppose that all vertices on $f$ except $v_1$, $v_5$ are bad. By symmetry, assume $v_2v_3$ is an edge in some 3-face, then $v_4v_5$, $v_5v_6$, $v_7v_8$ must be on 3-faces accounting for an $M$-face, contrary to Lemma 2.6.

Con 2. If the three willing vertices are not adjacent to each other, the bad vertices in $f$ can be grouped into two types, that is, group 1 is 1-1-3 and group 2 is 1-2-2. Otherwise, the bad vertices in $f$ can also be divided into two groups, that is, group 3 is 1-4 and group 4 is 2-3. Do the analysis as in Con 1 and Subcase 3.2, the four grouping cases all contradict the conditions of Con 2.

In conclusion, $ch^{\ast }(f)\geq 0$ when $f\in F(G)\setminus f_0$.$\hfill\blacksquare$

\begin{lemma}
	$ch^{\ast}(f_0)>0$.
\end{lemma}

\noindent\textbf{Proof. }Incontrovertibly, $ch^{\ast}( f_0)\geq ch_1(f_0)$. Let $n_{Fa}(G)$ denote the number of $FA$-structures in $G$, then $0\leq n_{Fa}(G)\leq 3$ is much in evidence.

\textbf{Case 1.} Suppose that $n_{Fa}(G)>0$, the case where $f_0$ gives the most charge is that when $n_{Fa}(G)=1$ and all the other vertices on $C_0$, except for this $FA$-structure, have degree 2. By R4, $f_0$ sends $\frac{4}{3}$ to $d(f_0)-2$ 2-vertices and special vertices, and 1 to other two vertices on $C_0$, which implies $ch_1(f_0)\geq d(f_0)+4-\frac{4}{3}\times (d(f_0)-2)-1\times 2=\frac{1}{3}(14-d(f_0))\geq 0$. Besides, $ch_1(f_0)=0$ only if $f_0$ meets the above requirements along with $d(f_0)=14$. But $G$ has one $10^{+}$-face $f_1$ containing nine 2-vertices with $ch_1(f_1)>0$ as Lemma 3.2. By R5, $ch^{\ast }(f_0)= ch_1(f_0)+ch_1(f_1)>0$.

\textbf{Case 2.} If $n_{Fa}(G)=0$, then $f_0$ sends $\frac{4}{3}$ only to 2-vertex by R4. Since $G\ncong G[V(C_0)]$ and $G$ is 2-connected, $C_0$ has at least two $3^{+}$-vertices. Thus $f_0$ gives at most $\frac{4}{3}$ to each of the $d(f_0)-2$ 2-vertices and 1 to two $3^{+}$-vertices, and then $ch_1( f_0)\geq d(f_0)+4-\frac{4}{3}\times (d(f_0)-2)-1\times 2=\frac{1}{3}(14-d(f_0))$. $ch_1(f_0)=0$ only when $d(f_0)=14$ and $C_0$ accurately contains twelve 2-vertices and two $3^{+}$-vertices. But there is at least one $10^{+}$-face $f_2$ possessing 2-vertices in $G$, and $ch_1(f_2)>0$ in terms of Lemma 3.2. By R5, we have $ch^{\ast}(f_0) \geq ch_1(f_0)+ch_{1}(f_2)> 0 $.

According to the above analysis, $ch^{\ast }(f_0)>0$.$\hfill\blacksquare$

\end{document}